%Paper: alg-geom/9306009
%From: enrique@sunal1.mat.ucm.es (Enrique Arrondo)
%Date: Wed, 23 Jun 93 11:24:01 +0200

%This is a plain TeX document, typeset by Textures 1.5, 13 pages.
\magnification 1200
\baselineskip 14pt

\def\mapright#1{\smash{
   \mathop{\longrightarrow}\limits\sp {#1}}}

\def\mapdown#1{\Big\downarrow
   \rlap{$\vcenter{\hbox{$\scriptstyle#1$}}$}}

\centerline{CLASSIFICATION OF SMOOTH CONGRUENCES WITH A FUNDAMENTAL CURVE}

\centerline{by {\it E. Arrondo}\footnote*{Departamento de Algebra,
Facultad de Ciencias Matem\'aticas, Universidad Complu\-tense de Madrid,
28040 Madrid, Spain}, {\it M. Bertolini}\footnote{**}{Dipartimento di
Matematica, Universit\`a degli Studi di Milano, via C. Saldini, 50,
20133 Milano, Italy} and {\it C. Turrini**}}
\bigskip
\bigskip
{\bf Abstract:} We give a classification and a construction of all smooth
$(n-1)$-dimen\-sional varieties of lines in ${\bf P}\sp n$ verifying that all
their lines meet a curve. This also gives a complete classification of
$(n-1)$-scrolls over a curve contained in $G(1,n)$.
\bigskip
\bigskip

{\bf Introduction.}

The geometry of the line, i.e. the study of the line as the main element,
was a very popular subject at the end of last century and beginning of
this. A particular attention was given to varieties of lines, also called
{\it congruences}. In other words, a congruence is a subvariety of a
Grassmannian $Gr(1,{\bf P}\sp n)$ of lines in ${\bf P}\sp n$. When we will
refer
to a congruence, we will restrict ourselves to the case of subvarieties of
dimension $n-1$ (hence equal to the codimension).

  A lot of papers were published at that time by several mathematicians
about congruences of lines in ${\bf P}\sp 3$. Of important relevance is the
work of G. Fano, who made several classifications of these congruences
under different assumptions. One fixed assumption in his classification
was the non-existence of what was called a {\it fundamental curve}.
A fundamental curve is a curve in ${\bf P}\sp 3$ such that all lines of
the congruence meet it. The reason for excluding this possibility does
not seem to be clear.
Recently, after congruences became popular again, M. Gross and the first
author classified all those congruences excluded by Fano, i.e they gave a
classification of all smooth congruences of lines in ${\bf P}\sp 3$ with a
fundamental curve.

  Not only the congruences of lines in ${\bf P}\sp 3$, but also in ${\bf
P}\sp 4$ (i.e. threefolds in $Gr(1,{\bf P}\sp 4)$) came back to the attention
of nowadays mathematicians. When the last two authors started the
classification of smooth congruences of lines in ${\bf P}\sp 4$ with a curve
meeting all lines, they found that the same key numerical relation as for
${\bf P}\sp 3$ holds (see lemma 4). This was the starting point of this joint
work dealing with the general case of smooth congruences of lines in any
${\bf P}\sp n$ with all lines meeting a curve (in the sequel such a curve will
be called fundamental, see section 0 for the definition).

  In section 1, we deal with the case in which the curve is a line. This
is a particular case in the sense that there are infinitely many families
of these congruences, while for other curves we find, for each $n$, a
finite number of families. On the other hand, this is a generalization of
the result in case $n=3$, which is known under another formulation.
Indeed, the Grassmannian $Gr(1,{\bf P}\sp 3)$ can be considered, under the
Pl\"ucker embedding, as a smooth quadric in ${\bf P}\sp 5$. Since the variety
of the lines meeting a given line is a singular hyperplane section,
congruences with a fundamental line can be viewed as surfaces in ${\bf
P}\sp 4$ contained in a singular quadric. It is then a classical result that
goes back to Roth (see [R] \S3) that such a smooth surface is either the
complete intersection of the quadric and another hypersurface (if it
misses the singular point of the quadric) or the rest of a plane under
such a complete intersection. We get an analogous result for any $n$.

  In section 2, we give the list of all possible smooth congruences. For the
proof we needed to refine some arguments of [A-G]. The base idea is to
lift the congruence as a hypersurface in
a desingularization of the variety of lines meeting a curve. Then one
obtains a numerical relation from the double-point formula (coming from
the fact that we have a smooth subvariety of $Gr(1,{\bf P}\sp n)$ with same
dimension as codimension), which together with Castelnuovo bound for the
genus in two different contexts allows us to conclude.

  In section 3 we give an explicit construction for all possible
congruences we found in section 2, which completes the classification.
We point out that all examples for $n=3$ (except the congruences of
bisecants) extend to any dimension.
  Moreover an easy argument shows that by the above constructions and [A]
we also obtain the complete classification of $(n-1)$-dimensional scrolls over
a curve, contained in $G(1,n)$.

   The results of this paper are just one way to generalize the classical
results on congruences with a fundamental curve.
Another natural and interesting generalization, which some of the
authors are now studying, are the congruences of lines in  ${\bf P}\sp n$
with higher dimensional fundamental locus, e.g. a surface with infinitely many
lines through each point of its.

   The first author wants to thank a Del Amo grant from Universidad
Complutense, which allowed him to stay at MSRI for the year 1992/93. The
second author also wants to thank G.N.S.A.G.A. of the C.N.R., which
allowed her visiting MSRI. They found there a wonderful work atmosphere
where they could start collaborating and also
completing the paper. We would also like to thank Mark Gross for several useful
conversations and ideas.

\bigskip

{\bf \S0. Preliminaries and notations.}

As said in the introduction, by a {\it congruence} we will mean an
$(n-1)$-dimensional variety (over the
complex numbers) in the
Grassmannian variety $G(1,{\bf P}\sp n)$ of lines in ${\bf P}\sp n$. We will
denote usually by $G$ to this Grassmannian.

We will write

\itemitem{$\Omega(A,B)=$} Schubert variety of lines in $G$ meeting $A$ and
contained in $B$, where $A$ and $B$ are linear subespaces of ${\bf P}\sp n$
such that $A\subset B$

\itemitem{$\Omega(i,j)=$} class, in the Chow ring of $G$,
of $\Omega(A,B)$, where dim$A=i$ and dim$B=j$

We will say that a point $P$ of ${\bf P}\sp n$ is a $k$-{\it fundamental point}
of a congruence $Y$ if there is a $k$-dimensional family of lines of $Y$
through
$P$.
A curve $C$ of ${\bf P}\sp n$ is a $k$-{\it fundamental curve} for $Y$ if all
its
points are $k$-fundamental points. In other words an $(n-2)$-fundamental curve
is a curve which is met by all the lines of the congruence; in the sequel we
will briefly call such a curve a {\it fundamental curve}.

Given a congruence $Y$, we will denote

\itemitem{$a=$} number of lines of $Y$ passing through a general
point of ${\bf P}\sp n$

\itemitem{$b=$} number of lines of $Y$ contained in a general hyperplane $H$
and meeting a general line of $H$

\itemitem{$d=$} degree of the curve $C$

\itemitem{$g=$} geometric genus of $C$

\itemitem{$e=$} degree of the cone formed by all lines
of $Y$ passing through a general point of $C$

We will not consider the case when $Y$ consists of all lines passing
through one point, so that we will assume $e\ge1$ for any fundamental
curve. Note that, although the
$(n-1)$-dimensional cycles of $G$ are generated by $\{\Omega(i,n-i)\mid
i=0,\ldots ,[{{n-1}\over 2}]\}$, in the case of congruences with
fundamental curve the only non-zero intersections with the
class of $Y$ are $a=[Y]\cdot \Omega(0,n)$ and $b=[Y]\cdot \Omega(1,n-1)$.
\bigskip

{\bf \S1. Smooth congruences with a fundamental line.}

For this section, we will view $G$ as embedded in a projective space under
the Pl\"ucker embedding. In this framework, a hypersurface of degree $l$
will mean the intersection of $G$ with such a hypersurface in the Pl\"ucker
embedding.

\proclaim Theorem 1. Let $Y\subseteq G(1,n)$ be a smooth congruence such
that all its lines meet a given line $\Lambda\subseteq{\bf P}\sp n$.
Denote by $\Gamma$ the cone $\Omega(\Lambda,{\bf P}\sp n)$.
Then either
\item{j)} $b=a$ and $Y$ is the intersection of  $\Gamma$ with a hypersurface
of degree $a$\hfill
or
\item{jj)} $b=a-1$ and $Y$ is linked to an $(n-1)-$fold of degree $n-2$
passing through
the vertex of $\Gamma$, in the intersection of $\Gamma$ with a hypersurface
of degree $a+1$\hfill or
\item{jjj)} The dimension is $n=3$, $a=b-1$ and $Y$ is linked under the
complete intersection of $\Gamma$ and a hypersurface of degree $b$ to
a plane consisting of all lines passing through a point.

\vskip 5truemm
{\it Proof:} Since the case $n=3$ is already known, we will assume $n\ge4$.
A modern proof for the case $n=3$ can be found in [G1], from where we took
the idea for this general proof.

It is easy to see that the Schubert variety
${\it\Gamma}=\Omega(\Lambda,{\bf P}\sp n)$ is a cone over the $(n-1)$-fold
$B={\bf P}\sp 1\times{\bf P}\sp
{n-2}$
with vertex the point $\lambda\in G(1,n)$ corresponding to ${\it\Lambda}$.
Consider
the blow up ${\it\Gamma\sp *}$ of ${\it\Gamma}$ at $\lambda$; then
${\it\Gamma\sp *}
\simeq{\bf P}$
($O_B \oplus O_B(-1)$) and the exceptional divisor $E$ is isomorphic to $B$.
Put ${\cal E} = O_B \oplus O_B(-1)$ and $p:{\cal E}\rightarrow B$;
Then, with standard notations, $Pic(B)$ ($Pic (E)$ resp.)
is generated by ${\bf P}\sp 1_1\times{\bf P}\sp {n-3}_2$ and ${\bf P}\sp
0_1\times{\bf
P}\sp {n-2}_2$ ( ${\bf P'}\sp 1_1\times{\bf P'}\sp {n-3}_2$ and ${\bf P'}\sp
0_1\times{\bf
P'}\sp {n-2}_2$  resp.), so that $Pic({\it\Gamma\sp *})$ is generated by $E$,
$Q=p*({\bf P}\sp 1_1\times{\bf P}\sp {n-3}_2)$ and
$Z={\bf P}\sp 0_1\times{\bf P}\sp {n-2}_2$.
Now, for a $(n-1)$-fold $Y\sp *\subseteq {\it\Gamma\sp *}$ to be mapped to a
non
singular
$(n-1)$-fold $Y\subseteq {\it\Gamma}$ there are two possibilities:

i) $Y\sp *$ does not intersect $E$, and then the blow up map is an isomorphism
between $Y\sp *$ and $Y$;

or

ii) $Y\sp *$ intersects $E$ and then $W = Y\sp *\cap E$ is contracted to the
vertex of ${\it\Gamma}$ which must be a smooth point of $Y$. So W must
be a divisor in $Y\sp *$ isomorphic to ${\bf P}\sp {n-2}$.

\noindent Now suppose that, in $Pic({\it\Gamma\sp *})$, it is $Y\sp * = \alpha
E +
 \beta Q
 + \gamma Z$, so that the numerical equivalency class of $W$ in $E$ is

\centerline{$W = (\beta -\alpha)({\bf P'}\sp 1_1\times{\bf P'}\sp {n-3}_2) +
(\gamma -
\alpha)({\bf P'}\sp 0_1\times{\bf P'}\sp {n-2}_2).$}

In case $j)$, one has $\alpha=\beta=\gamma$, so that $Y$ is the complete
intersection of ${\it\Gamma}$ with an hypersurface of degree $\alpha = a$
(case $i)$).
In case $jj)$, being $n-2 > 1$, the only possibility for a divisor $W$ to be a
${\bf P}\sp {n-2}$ embedded into $E = {\bf P}\sp 1\times{\bf P}\sp {n-2}$ is
that
$W = {\bf P}\sp 0_1\times{\bf P}\sp {n-2}_2$, so that $\beta = \alpha$ and
$\gamma =
\alpha + 1$. Hence $Y\sp * = \alpha E +\alpha Q + (\alpha +1)Z$ and $Y$ is
linked,
in the complete intersection of ${\it\Gamma}$ with an hypersurface of degree
$\alpha + 1$, to the image of $E + Q$ in ${\it\Gamma}$, i.e. to an $(n-1)$-fold
of degree $n-2$ passing through the vertex of ${\it\Gamma}$ (case $ii)$).
\bigskip

{\bf \S2. Smooth conguences with a fundamental curve}

In this section we will prove the following

\proclaim Theorem 2. Let $Y$ be a smooth congruence having an integral
curve $C$ as a fundamental curve. Then one of the following holds:
\item{(i)} $n=3$ and the congruence consists of the bisecants to either a
twisted cubic or an elliptic quartic in ${\bf P}\sp 3$.
\item{(ii)} The curve $C$ is a line (and $Y$ is described in Theorem $1$).
\item{(iii)} The congruence is a scroll (i.e. $e=1$) and either $C$ is a
conic and $a=1$ or $2$, $b=2$, or $C$ is a plane cubic and $a=b=3$
\item{(iv)} The curve $C$ is a plane cubic, $e=2$, $a=3$ and $b=6$

We will follow several steps in the proof, that we state as lemmas.

\proclaim Lemma 3. The fundamental curve $C$ is smooth.

 {\it Proof:} Consider the subscheme $V$ of $G$
corresponding to lines in ${\bf P}\sp n$ meeting $C$. It is singular, and its
singular locus is the locus of bisecants to $C$. This locus has several
components: one for each singular point of $C$ (consisting of all lines
passing through that point) and another one of the expected dimension
two given by the closure of all bisecant lines at smooth points of $C$.
Hence $Y$ is not contained in the singular locus of $V$ unless dim $Y=2$
and $Y$ is the congruence of bisecants to $C$. This case has
already be studied in [A-G], and corresponds to (i) in the theorem.

Let $f:\tilde C\to C$ be the normalization of $C$, and
$L=f\sp *({\cal O}_C(1))$. The scheme $V$ has a natural
desingularization given by
$X={\bf P}(f\sp *({\Omega_{{\bf P}\sp n}}_{\mid C}(2))$.
Call $p:X\to\tilde C$ the natural projection. Since by definition of
fundamental curve $Y$ is contained in $V$, it has a  lift $\tilde Y$ to
$X$. If we are not in case (i), then the map $\pi:\tilde Y\to Y$ is
birational. If it is not an isomorphism, it must contract a curve $E$
of $\tilde Y$ (just applying the Zariski main theorem, because $Y$ is
smooth). The curve $E$ meets all fibers of $p_{\mid\tilde Y}$, since
none of these fibers is contracted by $\pi$. This easily implies that
the line represented by the image of $E$ under $\pi$
is $C$, which is case (ii) in the theorem.

Hence, assuming that $C$ is not a line we have that $\pi$ is in fact an
isomorphism between $\tilde Y$ and $Y$. This implies that $C$ is smooth.
Indeed, if two different points $p_1$ and
$p_2$ of $\tilde C$
go to the same point of $C$, then the images of the cones at $p_1$ and
$p_2$ must have some common line. This contradicts the fact that $\pi$ is an
isomorphism. The same argument holds when $p_1$ and $p_2$ are infinitely near
points. A precise proof can be found in [G2]. The idea is that the tangent
space
at a point of $\tilde Y$ (given by a point $p$ of $\tilde C$ and a line $l$
through it) splits as the sum of the tangent space of $\tilde C$ at $p$ and the
tangent space at $l$ of the cone with vertex $p$. Since this maps
isomorphically
to the tangent space of $Y$, hence the tangent space of $\tilde C$ at $p$ also
maps isomorphically to the tangent space of $C$ at $f(p)$.

\proclaim Lemma 4. Except for cases (i) and (ii), there is a numerical
relation
$$2d\sp 2e\sp 2-4de\sp 2-2e\sp 2g-2de-2eg+2e\sp 2+2e+D(1+2e-2de)+D\sp 2=0$$
(This is the same as (*) in [G2] p. 137).

{\it Proof:} The map from $X$ to $G$ is given as
follows.

First, look at the commutative diagram of exact sequences defining
$F$ as a push-out (we will write ${\cal O}_X(a,b)$ for
$p\sp *({\cal O}_C(a))\otimes{\cal O}_X(b)$):
$$\matrix{&&0&&0&&&&\cr
          &&\mapdown{}&&\mapdown{}&&&&\cr
          &&\Omega_{X/C}(-1,1)&=&\Omega_{X/C}(-1,1)&&&&\cr
          &&\mapdown{}&&\mapdown{}&&&&\cr
          0&\mapright{}&p*{\Omega_{{\bf P}\sp n}}_{\mid C}(1,0)&\mapright{}
&H\sp 0({\cal O}
_{{\bf P}\sp n}(1))\otimes{\cal O}_X&\mapright{}&{\cal O}_X(1,0)&\mapright{}
&0\cr
          &&\mapdown{}&&\mapdown{\alpha}&&\mid \mid&&\cr
          0&\mapright{}&{\cal O}_X(-1,1)&\mapright{}&F&\mapright{}&{\cal
O}_X(1,0)&\to&0\cr
          &&\mapdown{}&&\mapdown{}&&&&\cr
          &&0&&0&&&&\cr}\eqno{(1)}$$
   The surjection $\alpha$ gives the map from $X$ to $G$; in particular
$\alpha$ is the pullback on $X$ of the canonical map from
$H\sp 0({\cal O}_{{\bf P}\sp n}(1))\otimes {\cal O}_G$ to $\cal Q$,
where $\cal Q$ is the universal bundle on $G$.
The Chow ring of $X$ is generated by $Pic(C)$ and $t=c_1({\cal
O}_X(1))$, which canonically identifies with the pullback of the hyperplane
section of $G$. The relation among them is given by $t\sp n=(n-1)Lt\sp {n-1}$,
and the class of a point is $Lt\sp {n-1}$. Hence,
the class of $Y$ in $X$ can be written as $et-D$, where $D$ is,
by abusing notation, the pullback of a certain divisor $D$
of $C$ (also, when no confusion arises, we will use $D$ for denoting the
degree of the divisor). Let us compute from these data the double-point
formula for $Y$. Writing $N=N_{Y/G}$, we have that
$a\sp 2+b\sp 2={\rm deg}(c_{n-1}(N))$. Let us compute each of
the members in the formula (most relations will come from the exact
sequences in diagram (1)).

The pullback to $X$ of the set of lines passing through a point of ${\bf
P}\sp 3$ is the locus where $n$ sections of the bundle $F$ are dependent.
We can apply then Porteous formula (see, for example [F] Thm. 14.4, from
where we also keep the notation for the Schur polynomial)
and get that
$$a=\Delta_1\sp {(n-1)}(F_{\mid Y})\cap[Y]=(t\sp {n-1}-(n-2)Lt\sp {n-2})
(et-D)=ed-D$$ (where we
used the notations in [F], and the fact that in our case
$\Delta_1\sp i(F)=t\sp i-(i-1)Lt\sp {i-1}$, as can be easely checked by
induction
on $i$)

The second degree $b$ is easily computed by geometrical means. It is
the number of lines of $Y$ contained in a general hyperplane $H$ of
${\bf P}\sp {n-1}$
and meeting a general line $r$ of $H$. The intersection of $H$ with $C$
consists of $d$ points, and through each of them, lines of $Y$ form a
cone of degree $e$, thus meeting $L$ in $e$ points. Hence $b=ed$.

Let us now compute the total Chern polynomial of $N$,
$c(N)=c({T_G}_{\mid Y})c(T_Y)\sp {-1}$. We use the fact that $T_G=F\otimes
T_{X/C}(1,-1)$ and the fact that $F$ appears as an extension in (1), as
well as the fact that $N_{Y/X}={\cal O}_X(et-D)$ to conclude that
$$c(N)=c(p\sp *({T_{{\bf P}\sp n}}_{\mid C}))c({\cal O}_X(2,-1))\sp {-1}
c({\cal O}_X(et-D))c(p\sp *T_C)\sp {-1}.$$
Now, it is a straightforward calculation to show that
$$c_{n-1}(N)=(e+1)t\sp {n-1}+(5e-ne-n+3)Lt\sp {n-2}+(e+1)Kt\sp {n-2}-Dt\sp
{n-2}$$
(where $K$ represents the canonical divisor of $C$). Now, its degree comes
from multiplying in the Chow ring of $X$ the above class with the class
$et-D$ of $Y$, to obtain
$${\rm deg}(c_{n-1}(N))=(4e\sp 2+2e)d+e(e+1)(2g-2)-(2e+1)D$$
Identifying this with $a\sp 2+b\sp 2$ and substituting the values of $a$ and
$b$ we finally get
$$2d\sp 2e\sp 2-4de\sp 2-2e\sp 2g-2de-2eg+2e\sp 2+2e+D(1+2e-2de)+D\sp 2=0$$

\proclaim Lemma 5. The curve $C$ is plane and hence the relation in
lemma 4 becomes
$$de(d-1)(e-1)=D(2de-D-2e-1) \eqno(2)$$

{\it Proof:} It is the same as in [G2] by getting a contradiction between
the relation in lemma 4 and Castelnuovo bound for the genus of non-plane
curves.

\proclaim Lemma 6. In the same hypothesis of lemmas 4 and 5, $b\le 2a$
and hence $e\le 3$.

{\it Proof:} This is just as in [A-G], by applying Hurwitz theorem to the
map from the curve of lines of $Y$ in $P$ (the plane containing $C$) to
$C$. The second inequality is a consequence of the first together with (2).
\vskip 5 truemm
We now complete the proof of the theorem by analyizing separately each of
the possible values for $e$.

1) If $e=1$, then the congruence $Y$ is a scroll of ${\bf P}\sp {n-2}$'s
over the curve $C$. More precisely, through each point $Q$ of $C$, the
lines of the congruence $Y$ are those contained in a given hyperplane
$H_Q$ of ${\bf P}\sp n$ containing $Q$. This distribution of hyperplanes
is given by the epimorphism in the exact sequence (defining the vector
bundle $E$):

$$0\to {\cal O}_C(D)\to\Omega_{{\bf P}\sp n}(2)\otimes{\cal O}_C\to
E\to 0\eqno(3)$$
(recall that $Y$ is the zero locus of
$p\sp *({\cal O}_C(D))\to p\sp *(\Omega_{{\bf P}\sp n}(2)\otimes{\cal O}_C)
\to{\cal O}_X(0,1)$).
Since the curve $C$ is plane, we can build the
following commutative diagram (defining the sheaf ${\cal F})$:
$$\matrix{&&&&0&&0&&\cr
  &&&&\mapdown{} &&\mapdown{} &&\cr
  &&&&{\cal O}_C(L)\sp {\oplus n-2}&=&{\cal O}_C(L)\sp {\oplus n-2}&&\cr
  &&&&\mapdown{} &&\mapdown{} &&\cr
  0&\to&{\cal O}_C(D)&\to&\Omega_{{\bf P}\sp n}(2)\otimes{\cal O}_C&\to&E&\to
&0\cr
  &&\mid\mid&&\mapdown{} &&\mapdown{} &&\cr
  0&\to&{\cal O}_C(D)&\to&\Omega_{{\bf P}\sp 2}(2)\otimes{\cal
O}_C&\to&{\cal F}&\to&0\cr
  &&&&\mapdown{} &&\mapdown{} &&\cr
  &&&&0&&0&&} \eqno(4)$$
The injectivity in the right column in (4) comes because otherwise the sheaf
${\cal F}$ would have rank two, just implying that all lines in
the plane $P$ are in the congruence.

Putting $e=1$ in expression
(2), we get that either deg$D=0$ (and hence $a=b=d$) or $a+b=3$. In this
second case, from
lemma 6 the only possiblity is $a=1$, $b=d=2$. So we assume that
deg$D=0$.

Assume the sheaf ${\cal F}$ appearing in (4) has torsion at some point $p$
of $C$.
This would imply that for each other $q$ in $C$, the line joining $p$ and $q$
is in the congruence. For $d>2$ this contradicts the fact proved in lemma 3
that $\pi$ is an isomorphism. Hence, ${\cal F}$ is an invertible sheaf in this
case, so it must be ${\cal F}={\cal O}_C(L-D)$, by looking at the first
Chern class. This gives the immersion
of $C$ into $\check{\bf P}\sp 2$ defined by assigning to each point in $C$ the
only
line through it contained in $P$. For $d>3$ there is only one $g_d\sp 2$ on $C$
(see [A-C-G-H] page 56), so that $L=L-D$ and thus $D=0$.
This would imply that all lines in $P$ pass through a fix point of $C$, which
is
absurd. Hence, $d\le3$.

2) If $e=2$, it happens as in [A-G] that the curve (of degree $a$) $W$ in the
projective space of hyperquadrics in ${\bf P}\sp n$ must be plane (just using
Castelnuovo bound together with lemma 5). This plane
cannot be contained in the locus $K$ of singular quadrics whose vertex is
at a given plane $P$ (same proof as in [A-G], using that the embedded
tangent space of $K$ at a cone contains the space of quadrics passing
through the vertex of the cone). This locus $K$ is defined by
the maximal minors of a $3\times (n+1)$ matrix of linear forms. Hence,
$a\le 3$, and the only numerical solutions for $(a,b)$ are $(2,6)$ and
$(3,6)$. The first one is ruled out by lemma 6.

3) If $e=3$, the curve $W$ in the space of hypercubics is contained in a
linear space $A$ of dimension three. This has to meet $K$ (space of cones
with vertex at the plane $P$) along a curve. Indeed, the embedded tangent
space to $K$ at a cubic cone is contained in the space of cubic
hypersurfaces that are singular at the vertex of the cone. Assume that
the intersection of $K$ and $A$ contains a surface $S$. Take a point
$p$ of $C$ and let $p'$ the point of $W$ representing the cone with
vertex $p$. By assuption, there is a plane in $A$ through $p'$ (the
tangent plane to $S$ at $p'$) consisting of singular cubics at $p$.
If $a>3$, this plane must meet $W$ in at least another point $q'$
(for general $p$). Hence the cone at the corresponding point $q$ of
$C$ is singular at $p$. When we vary the point $p$, if the point $q$
also varies, this means that all cones are singular, which is absurd
(the general fiber of the map $\tilde Y\to C$ is smooth). This means
that there is a point $q$ in $C$ whose cone through it is singular
along all points of $C$. In particular the pencil of
lines through $q$ and contained in $P$ should be in the congruence, which
is also absurd.

Now $K$ is defined by the
maximal minors of a $3\times{n+2\choose 2}$ matrix of linear forms (this
matrix corresponds to the derivatives of forms of degree 3 with respect to
the directions given by the plane $P$).
Hence, $W$ is contained in the locus of a ${\bf P}\sp 3$ given by the
maximal minors of a $3\times4$ matrix of linear forms, hence $a\le 6$.
The only numerical solutions of (2) are then
$(a,b)=(3,12),(4,12)$, which are again impossible by lemma 6.
\bigskip

{\bf \S3. Explicit constructions.}

\bigskip

\proclaim Lemma 7. Let $C\subseteq{\bf P}\sp n$ be a non singular plane curve
of degree $d$. Then
$\Omega_{{\bf P}\sp n}(2)_{|C}\cong(\oplus \sp {n-2}{\cal O}_C(L))\oplus{\cal
S}(P)$
where $L$ is an hyperplane section of $C$,  $E$ is a divisor on $C$
of degree 1 and ${\cal S}$ is a normalized rank 2 vector bundle of degree $d-2$
on $C$ with $\wedge\sp 2{\cal S}\cong{\cal O}_C(L-2P)$.
\vskip 1truecm
{\it Proof:}
{\rm First of all, notice that}
\vskip 5truemm

\centerline{$\Omega_{{\bf P}\sp n|{\bf P}\sp 2}\cong
(\oplus \sp {n-2}{\cal O}_{{\bf P}\sp 2}(-1))\oplus\Omega_{{\bf P}\sp 2}$.}
\vskip 5truemm

Twisting by ${\cal O}(2)$ and restricting to $C$, we get
\vskip 5truemm

\centerline{$\Omega_{{\bf P}\sp n}(2)_{|C}\cong(\oplus \sp {n-2}{\cal O}_C(L))
\oplus\Omega_{{\bf P}\sp 2}(2)_{|C}$.}
\vskip 5truemm
If $\Lambda$ denotes the plane containing $C$ and
\vskip 5truemm

\centerline{$f:X=\Omega_{{\bf P}\sp n}(2)_{|C} \to G(1,n)$}
\vskip 5truemm

\noindent is the projection, denote by $Y_0=f\sp {-1}(\Omega(1,\Lambda))$
the ruled surface of
lines of $X$ contained in $\Lambda$. Let ${\cal S}$ be the normalized
bundle defined
by $Y_0={\bf P}({\cal S})$. By the above, $Y_0$ can be identified with
${\bf P}(\Omega_{\Lambda}(2)_{|C})\subseteq {\bf P}(\Omega_{{\bf P}\sp n}(2)_
{|C})$, which implies $(3)$.
Now, with the same arguments as in [A-G], one shows that
$deg{\cal S}=d-2$, and that the numerical equivalence class of the divisor
inducing the map $Y_0\to\Omega(1,\Lambda)$ is $C_0+f$. By the way notice
also that the class of $Y_0$ in $X$ is the $(n-2)-th$ Chern class of
$(\oplus \sp {n-2}{\cal O}_X(t-p\sp *L))$.
\bigskip

{\bf Construction of case $(3,3)$}
\vskip 3truemm

Here $\Omega_{{\bf P}\sp n}(2)_{|C}\cong(\oplus \sp {n-2}{\cal
O}_C(L))\oplus{\cal F}
(P)$ as in Lemma 7, with $deg L = 3$, $deg P = 1$, $deg {\cal F}=1$ and
$\wedge\sp 2{\cal F}\cong{\cal O}_C(L-2P)$, so that $Y_0$ is the elliptic ruled
surface of invariant $-1$

Notice that, $C$ being a plane curve, the image of $X$ in the grassmannian

$G(1,n) \in {\bf P}\sp {{n(n+1)\over 2}-1}$ is in fact contained in a
${\bf P}\sp {3n-4}$, so that, since $dim |t| = 3n-4$, the embedding
of $X$ into $G(1,n)$ is given by the whole linear system $|t|$.

As shown in Lemma $4$, the class of $Y$ in
$X={\bf P}(\Omega_{{\bf P}\sp n}(2)_{|C})$, can be written as $t-p\sp *D$, $D$
being a degree $0$ divisor on $C$.

As $D$ is a divisor of degree $0$, then $\Omega_{{\bf P}\sp n}(2)_{|C}(-D)$
is generated by global sections, so $|t-p\sp *D|$ is base point free and
a general member is smooth.

Recall from \S2 that we have an exact sequence given by the right column of
(4) and that we showed that ${\cal F}={\cal O}_C(L-D)$ and $D\neq 0$. This
implies that $Y={\bf P}(E)={\bf P}({\cal O}_C(L-D)\oplus
{\cal O}_C(L)\sp {n-2})$. Since it is mapped to the Pl\"ucker embedding of
$G$ by the complete very ample linear series $\mid t\mid$, then it is
an embedding in $G$.
\vskip 3truemm

For the geometric construction, we observe first that the hyperplanes formed
at each
point must all contain a fixed linear space $\Lambda$ of codimension
three, which is disjoint with the plane $P$ containing $C$. Indeed, this
comes from the splitting of $E$. Hence these hyperplanes are defined as
the span of $\Lambda$ and some a line in $P$. The way of constructing these
lines is given in [Go] in the following way. Take a group structure on $C$
such that the origin is an inflection point. Now take a point $\sigma$
such that $3\sigma=L-D$ and fix a point $R$ outside $C$. For each point
$p$ in $C$, let $A$ and $B$ the points of $C$ in the line defined by $R$
and $-p+2\sigma$. Then, the line associated to $p$ is the one defined by
$-A-\sigma$ and $-B-\sigma$. We do not know how to prove this directly, but
just as in [G1] checking that this construction corresponds to the above
description and that, for fixed $D$, the dimension is the same (dimension
two for the choice $R$ and also for a section of $\Omega(2)_{|C}(-D)$).
\bigskip

{\bf Construction of case $(2,2)$}
\vskip 3truemm

Let $C$ be a conic, contained in a plane $L$ in ${\bf P}\sp n$. By Lemma $7$,
the sheaf $\Omega_{{\bf P}\sp n}(2)_{|C}$ is isomorphic to
$\oplus \sp {n-2}{\cal O}_{{\bf P}\sp 1}(2)\oplus(\oplus \sp 2{\cal O}_{{\bf
P}\sp 1}(1))$.

As $X={\bf P}(\Omega_{{\bf P}\sp n}(2)_{|C})$ we can write
$X={\bf P}(\oplus \sp {n-2}{\cal O}_{{\bf P}\sp 1}(2)\oplus (\oplus\sp 2{\cal
O}_
{{\bf P}\sp 1}(1))$.
Let us consider now the linear system $|t|$. As $h\sp 0({\cal O}_X(t))=3n-2$
and $G(1,n)\subset{\bf P}\sp {{n(n+1)\over 2}-2}$, the map $f:X \to G(1,n)$
is induced by a $(3n-3)$-dimensional subspace $V$ of $H\sp 0({\cal O}_X(t)$.
More precisely, $V=H\sp 0({\cal O}_{{\bf P}\sp 1}(2)\sp {n-2})\oplus V'$, where
$V'$ is a three-dimensional subspace of $H\sp 0({\cal O}_{{\bf P}\sp 1}(1)\sp
2)$.
\vskip 3truemm

\proclaim Claim. A smooth scroll with $a=b=2$ comes from an element
$Y\in |t|$ of one of the following types:
\item{(i)}$Y={\bf P}(E)$ with $E=\oplus \sp {n-1}{\cal O}_{{\bf P}\sp 1}(2)$
i.e. ${\bf P}\sp 1 \times {\bf P}\sp {n-2}$ embedded
in ${\bf P}\sp {3n-4}$ via ${\cal O}(2,1)$
\item{(ii)}$Y={\bf P}(E)$ with
$E= \oplus \sp {n-3}{\cal O}_{{\bf P}\sp 1}(2)\oplus {\cal O}_{{\bf P}\sp 1}(1)
\oplus {\cal O}_{{\bf P}\sp 1}(3)$.
\vskip 3truemm

{\it Proof of claim}: The exact sequence (3) in \S2 becomes
\vskip 3truemm

$$0 \to {\cal O}_{{\bf P}\sp 1} \to \oplus \sp {n-2}{\cal O}_{{\bf P}\sp 1}(2)
\oplus(\oplus \sp 2{\cal O}_{{\bf P}\sp 1}(1)) \to E \to 0.$$
\vskip 3truemm
We look at the different possibilities for $E$ depending on the
composed map
$$\varphi:{\cal O}_{{\bf P}\sp 1} \to {\cal O}_{{\bf P}\sp 1}(2)\sp {n-2}
\oplus{\cal O}_{{\bf P}\sp 1}(1)\sp 2\to{\cal O}_{{\bf P}\sp 1}(1)\sp 2$$
\item{(i)} If $\varphi$ is an injective morphism of bundles, then $E$ is
an extension of coker$\varphi={\cal O}_{{\bf P}\sp 1}(2)$ by
${\cal O}_{{\bf P}\sp 1}(2)\sp {n-2}$, hence
$E= \oplus \sp {n-1}{\cal O}_{{\bf P}\sp 1}(2)$

\item{(ii)} If $\varphi$ is (after changing basis) zero on one factor
${\cal O}_{{\bf P}\sp 1}(1)$ and different from zero on the other, then $E$
contains the first ${\cal O}_{{\bf P}\sp 1}(1)$ as a direct summand.
The complement of this will be an extension of ${\cal O}_{{\bf P}\sp 1}(3)$
(this is the cokernel for a general monomorphism of bundles
${\cal O}_{{\bf P}\sp 1}$to${\cal O}_{{\bf P}\sp 1}(1)\oplus
{\cal O}_{{\bf P}\sp 1}(2)$ by ${\cal O}_{{\bf P}\sp 1}(2)\sp {n-3}$, so that
$E= \oplus \sp {n-3}{\cal O}_{{\bf P}\sp 1}(2)\oplus {\cal O}_{{\bf P}\sp 1}
(1))\oplus {\cal O}_{{\bf P}\sp 1}(3))$.

\item{(iii)} If $\varphi=0$, then $Y={\bf P}(E)$ must contain
${\bf P}(\Omega_{{\bf P}\sp 2}(2)_{|C})$, which would imply that all
the lines in $L$  through any point of $C$ are in the congruence.
This is a contradiction.
\bigskip
The case (ii) occurs on an open subset of the locus defined by a quadric
equation in $V'$. Hence, we can find a variety $Y$ in case (ii) not
defined by an element of $V$.
We want to prove now that if $Y$ is not in $V$, it is mapped
isomorphically in $G(1,n)$.

Let us consider the sequence
\vskip 3truemm

\centerline {$0 \to H\sp 0({\cal O}_X) \to H\sp 0({\cal O}_X(t))
\to H\sp 0({\cal O}_Y(t))\to 0$.}
\vskip 3truemm

Now $Y$ not being in $V$, $V$ is mapped onto $H\sp 0({\cal O}_Y(t))$ so that
the immersion of $Y$ into the grassmannian is determined by a complete
linear system, which is in fact very ample.

This concludes the proof.
\vskip 3truemm
We now furnish an example of geometric construction of the two possible
scrolls of bedegree (2,2) quoted in cases (i) and (ii) of the claim
above, generalizing some examples of [Go] to the case of higher
dimension.

Case (i). Let $L$ be a $2$-plane into ${\bf P}\sp n$ and $C$
a smooth conic on it, and let $\Lambda$ be a ${\bf P}\sp {n-3}$ skew with
$L$. For each point $p \in C$ take the tangent line $t_p$ to $C$ at $p$,
consider the ${\bf P}\sp {n-1}=[\Lambda , t_p]$ spanned by $\Lambda$ and
$t_p$ and denote by $\Sigma _p$ the ${\bf P}\sp {n-2}$ of the lines through
$p$  contained in $[\Lambda , t_p]$. A model for $Y$ is generated by the
$\Sigma _p$'s when $p$ varies in$C$.

Case (ii). Let $L$ and $C$ be as above and take a point
$p_0$ in $C$. Each ${\bf P}\sp {n-2}= \Sigma _p$ of the lines of the
congruence $Y$ through a point $p$ of $C$ generates a ${\bf P}\sp {n-1}$
passing through $p_0$, the ${\bf P}\sp {n-1}$ generated by $\Sigma _{p_0}$
is tangent to $C$ at $p_0$ but does not contain $L$.
\vskip 3truemm

{\bf Construction of case $(1,2)$}
\vskip 3truemm

In this case, we obtain that the exact sequence (3) in \S2 is
$$0\to{\cal O}_{{\bf P}\sp 1}(1)\to{\cal O}_{{\bf P}\sp 1}(1)\sp 2\oplus
{\cal O}_{{\bf P}\sp 1}(2)\sp {n-2}\to E\to 0$$
{}From this we immediately get (since the map ${\cal O}_{{\bf
P}\sp 1}(1)\to{\cal O}_{{\bf P}\sp 1}(1)\sp 2$ is not zero) that $Y={\bf P}(E)=
{\bf P}({\cal O}_{{\bf P}\sp 1}(1)\oplus {\cal O}_{{\bf P}\sp 1}(2)\sp {n-2})$.
The map to $G$ is given by the complete linear series of the canonical
${\cal O}(1)$ for $E$, which is very ample, so that we get a smooth
congruence.

The geometric interpretation is as follows. The curve in
$\check{\bf P}\sp n$ of the hyperplanes defined at points of $C$ has degree
$a=1$, so that all these planes must be contained in a linear subspace
$\Lambda$ of codimension two. Hence all lines in the congruence meet $C$
and $\Lambda$. The set of lines meeting both $C$ and $\Lambda$ defines a
congruence with $a=b=2$, so that it contains another congruence with $a=1$
and $b=0$. This implies that $\Lambda$ meets $C$ in one point and $Y$
consists of the closure of all lines meeting $\Lambda$ and $C$ in points
different from $C\cap\Lambda$.
\bigskip
{\bf Remark}: For $n\not= 4$ the only  $(n-1)$-dimensional scrolls over
a curve contained in $G(1,n)$ are the three ones constructed above (cases
$(3,3)$, $(2,2)$ and $(1,2)$). Indeed it is easy to see that, for $n\not= 4$,
the only linear subspaces of dimension $n-2$ contained in $G(1,n)$ are the
Schubert cycles $\Omega (0,n)$.

For $n=4$ the linear subspaces of dimension $2$ contained in $G(1,4)$ are the
Schubert cycles $\Omega (0,4)$ and $\Omega (1,3)$. Hence in this case one has
also to consider the scrolls whose fibers are $\Omega (1,3)$, described in [A],
which are without a fundamental curve and therefore not included in our paper.

This completes the classification.
\vskip 3truemm
{\bf Construction of case $(3,6)$}
\vskip 3truemm

Let $C\subseteq {\bf P}\sp n$ be the non singular plane cubic quoted in
theorem $2$, case $(iv)$. Then $e=2, d=3, D=3$. Now, with notation of
Lemma $7$ we have $deg{\cal F}=+1$, i.e. $Y_0$ is the elliptic ruled
surface of invariant $-1$.
As shown in Lemma $4$, the class of $Y$ in
$X={\bf P}(\Omega_{{\bf P}\sp n}(2)_{|C})$, can be written as $2t-p\sp *D$.

\noindent So we study the linear system ${\it L}=|2t-p\sp *D|$, on
$X$. We have
\vskip 5truemm

$h\sp 0({\cal O}_X(2t-p\sp *D))={3(n+1)(n-2)\over 2}+h\sp 0(S\sp 2({\cal
F})(2P-D)).$
\vskip 5truemm

Since $e=-1$, there exists a curve numerically equivalent to $2C_0-f$
(see [E] lemma 1.4),
where $C_0$ and $f$ denote respectively a fundamental section and a
fibre of $Y_0$, so that we can fix a $D$ such that
$h\sp 0({\cal O}_X(2t-p\sp *D)) \geq {3(n+1)(n-2)\over 2} +1$.

Now we show that the surface
$Y_0$ is not contained in the base locus of ${\it L}$. To do it, consider
the ideal ${\cal I}$ of $Y_0$ and the exact sequence

$$0 \to{\cal I}(2t-p\sp *D) \to {\cal O}_X(2t-p\sp *D) \to
 {\cal O}_{Y_0}(2t-p\sp *D) \to 0.$$

Notice that for $i>0$, $R\sp ip_*{\cal I}(2t-p\sp *D)=0$, since the fibre of
this sheaf at a point $p_0$ of $C$ is $H\sp i({\cal J}(2))$, ${\cal J}$
denoting the ideal of the ${\bf P}\sp 1$
defined by $Y_0$ in the fibre ${\bf P}\sp {n-1}$ of $X$ at $p_0$.
Therefore  one can apply $p_*$ to the exact sequence above to get the
sequence

$$0 \to p_*{\cal I}(2t-p\sp *D) \to S\sp 2(\Omega_{{\bf P}\sp n}(2)_{|C})(-D)
\to S\sp 2(\Omega_{{\bf P}\sp 2}(2)_{|C})(-D) \to 0$$

\noindent which splits; this easily implies that
$h\sp 0(p_*{\cal I}(2t-p\sp *D))=h\sp 0({\cal I}(2t-p\sp *D))=h\sp 0({\cal
O}_X(2t-p\sp *D))-1=
{3(n+1)(n-2)\over 2} $, so that $Y_0$ is not contained in the base locus
of ${\it L}=|2t-p\sp *D|$.

As $Y\in|2t-p\sp *D|$
and $[Y_0]=t\sp {n-2}-3(n-2)t\sp {n-1}p\sp *P$, as remarked in Lemma $7$
then $Y{\cap Y_0}$ is a curve $C'$
numerically equivalent to $2C_0-f$, since
$C'\sp 2=Y\sp 2{\cdot Y_0}= 4t\sp n-12(n-1)t\sp {n-1}p\sp *P =0.$

As $C'$ is the unique effective divisor in its linear equivalence class,
it is contained in the base locus of the system $|2t-p\sp *D|$. Now
restricting this system to a fiber  $p\sp {-1}(P)$, at a point $P{\in C}$,
we get the following exact sequence:

$$0 \to {\cal O}_X(2t-p\sp *(D+P)) \to {\cal O}_X(2t-p\sp *D)
\to {\cal O}_{p\sp {-1}(P)}(2) \to 0$$

\noindent by which we get that the image of $H\sp 0({\cal O}_X(2t-p\sp *D))$
in $H\sp 0({\cal O}_{p\sp {-1}(P)}(2))$ has dimension ${(n\sp 2+n-4)\over 2}$.
Hence $|2t-p\sp *D|$ cuts out on the fiber a system of dimension
${(n-2)(n+3)\over 2}$ of quadrics which has at least two base-points,
hence just the two points of $Y{\cap Y_0}{\cap p\sp {-1}(P)}$. Thus
$Y{\cap Y_0}$ is the base locus of $|2t-p\sp *D|$, and $Y$ is smooth
away from  $Y{\cap Y_0}$.

With the same argument as in [A-G] we can conclude that $Y$ is embedded in
$G(1,n)$ by the map $f:X{\to G(1,n)}$ as a smooth congruence of bidegree
$(3,6)$.

\bigskip

{\bf References:}

\item{[A]} Alzati, A., {\it 3-Scroll immersi in $G(1,4)$}, Ann. Univ.
Ferrara, {\bf 32} (1986), 45-54.
\item{[A-C-G-H]} Arbarello, E.-- Cornalba, M.-- Griffiths, P.--
Harris, J., {\it Geometry of Algebraic Curves}, Vol. I Grund. der Math.
Wissen. {\bf 267} Springer-Verlag (1985).
\item{[A-G]} Arrondo, E.-- Gross, M., {\it On smooth surfaces in
$Gr(1,{\bf P}\sp 3)$ with a fundamental curve}, to appear in Man. Math.
\item{[E]} Ein, L., {\it Non-degenerate surfaces of degree $n+3$ in
${\bf P}\sp n$}, Crelle Journal Reine Angew. Math. {\bf 351} (1984), 1-11.
\item{[F1]} Fano, G., {\it Sulle congruenze di rette del terzo ordine
prive di linea singolare}, Att. Acc. di Scienze Torino {\bf 29} (1984),
474-493.
\item{[F2]} Fano, G., {\it Nuove ricerche sulle congruenze di rette del
$3\sp 0$ ordine prive di linea singolare}, Memoria della Reale Acad. Sc.
Torino (2) {\bf 51} (1902), 1-79.
\item{[F]} Fulton, W., {\it Intersection Theory}, Ergebnisse (3) {\bf 2},
Springer (1984).
\item{[Go]} Goldstein, N. {\it Scroll surfaces in $Gr(1,{\bf P}\sp 3)$},
Conference on Algebraic Varieties of small dimension (Turin 1985), Rend.
Sem. Mat. Univ. Politecnica, Special Issue (1987) 69-75.
\item{[G1]} Gross, M., {\it Surfaces in the four-dimensional
Grassmannian}, Ph. D. thesis, Berkeley (1990).
\item{[G2]} Gross, M., {\it The distribution of bidegrees of smooth
surfaces in $Gr(1,{\bf P}\sp 3)$}, Math. Ann. {\bf 292} (1992), 127-147.
\item{[R]} Roth, L., {\it On the projective classification of surfaces},
Proc. London Math. Soc., {\bf 42} (1937), 142-170.

\bigskip
\centerline{E-mail of the authors:}

\centerline{arrondo@mat.ucm.es}

\centerline{bertolin@vmimat.mat.unimi.it}

\centerline{turrini@vmimat.mat.unimi.it}

\end